\documentclass[letterpaper, 10 pt, conference]{ieeeconf}  

\IEEEoverridecommandlockouts                              
\usepackage{amsfonts}
\usepackage{graphicx,xcolor}
\usepackage{amsmath}
\usepackage{amssymb}
\usepackage{dsfont}

\usepackage{ifthen}
\usepackage{color}
\usepackage{cite}

\usepackage{mathtools}
\usepackage{booktabs}
\usepackage{url}
\usepackage{hyperref}
\usepackage[font={small}]{caption}
\usepackage{float}
\usepackage{multirow}

\usepackage{caption}
\usepackage{subcaption}
\usepackage{comment}

\def\scr#1{{\cal #1}}
\newcommand{\R}{{\rm I\!R}}
\def\eq#1{\begin{equation}#1\end{equation}}

\newcommand{\bbb}{\mathbb}
\newtheorem{theorem}{Theorem}
\newtheorem{lemma}{Lemma}

\newtheorem{proposition}{Proposition}
\newtheorem{corollary}{Corollary}

\newcommand{\dfb}{\stackrel{\Delta}{=}}
\def\qed{ \rule{.08in}{.08in}}

\newcommand{\1}{\mathbf{1}}

\title{\LARGE \bf Fast Consensus over Almost Regular Directed Graphs
}

\author{Susie Lu \hspace{.3in} Marco Gamarra \hspace{.3in} Ji Liu
\thanks{
S. Lu is with Stanford Online High School (SusieLu@ohs.stanford.edu).
M. Gamarra is with the Air Force Research Laboratory
(marco.gamarra@us.af.mil).
J. Liu is with the Department of Electrical and Computer Engineering at Stony Brook University
(ji.liu@stonybrook.edu).
}
}

\begin{document}

\maketitle
\thispagestyle{empty}
\pagestyle{empty}


\begin{abstract}
This paper studies an open consensus network design problem: identifying the optimal simple directed graphs, given a fixed number of vertices and arcs, that maximize the second smallest real part of all Laplacian eigenvalues, referred to as algebraic connectivity. For sparse and dense graphs, the class of all optimal directed graphs that maximize algebraic connectivity is theoretically identified, leading to the fastest consensus. For general graphs, a computationally efficient sequence of almost regular directed graphs is proposed to achieve fast consensus, with algebraic connectivity close to the optimal value.
\end{abstract}


\section{Introduction}

In a network of $n$ agents, consensus means that all $n$ agents reach an agreement on a specific value for their agreement variables, $x_i(t)\in\R$, $i\in\{1,2,\ldots,n\}$. 
A continuous-time linear consensus process over a simple directed graph $\bbb G$ can be typically modeled by a linear differential equation of the form $\dot x(t)=-Lx(t)$, where $x(t)$ is a vector in $\R^n$ and $L$ is the ``Laplacian matrix'' of $\bbb G$ \cite{reza1}. 
It is well known that a consensus will be reached if and only if the underlying graph $\bbb G$ is ``rooted'' \cite[Theorem~2]{RenChapter}. 

For any simple directed graph $\bbb G$ with $n$ vertices, we use $D(\bbb G)$ and $A(\bbb G)$ to denote its in-degree matrix and adjacency matrix, respectively. Specifically, $D(\bbb G)$ is an $n\times n$ diagonal matrix whose $i$th diagonal entry equals the in-degree of vertex $i$, and $A(\bbb G)$ is an $n\times n$ matrix whose $ij$th entry equals 1 if $(j,i)$ is an arc (i.e., a directed edge) in $\bbb G$ and otherwise equals 0. The in-degree based Laplacian matrix of $\bbb G$ is denoted and defined as $L(\bbb G)=D(\bbb G)-A(\bbb G)$. There exist other definitions of the Laplacian matrix for a directed graph. For example, the out-degree based Laplacian matrix is defined by replacing the in-degree matrix $D$ with the out-degree matrix and the adjacency matrix $A$ with its transpose \cite{wiki}.
This paper focuses on in-degree based Laplacian matrices, which we refer to simply as Laplacian matrices.
It is easy to see that any Laplacian matrix has an eigenvalue at 0 with an eigenvector $\1$, where 
$\1$ is a column vector in $\R^n$ whose entries all equal 1. Using the Gershgorin circle theorem \cite{gershgorin}, it is straightforward to show that all Laplacian eigenvalues, except for those at 0, have positive real parts, as was done in \cite[Theorem 2]{reza1} for out-degree based Laplacian matrices. Therefore, the smallest real part of all Laplacian eigenvalues is always 0.

From standard linear systems, the continuous-time linear system $\dot x(t)=-Lx(t)$ will reach a consensus as $t\rightarrow\infty$ (i.e., $\lim_{t\rightarrow\infty} x(t) = a\1$ where $a$ is a constant) exponentially fast if and only if the second smallest real part among all eigenvalues of $L$ is positive, which determines the worst-case convergence rate of consensus.

In the special case when $\bbb G$ is a simple undirected graph, each undirected edge between two vertices $i$ and $j$ can be equivalently replaced by a pair of directed edges $(i,j)$ and $(j,i)$; then $L(\bbb G)$ is a symmetric matrix and thus has a real spectrum. Its second smallest eigenvalue is called the algebraic connectivity of $\bbb G$ and is positive if and only if $\bbb G$ is connected \cite{Fiedler73}. 
Motivated by this, we call the second smallest real part among all Laplacian eigenvalues of a directed graph $\bbb G$ as the algebraic connectivity of $\bbb G$ and denote it by $a(\bbb G)$.
This quantity has been studied in \cite{Asadi16,Asadi16cdc,Zhang17}.

From the preceding discussion, it is of particular interest to identify the optimal directed graphs, given a fixed number of vertices and arcs, that maximize algebraic connectivity, leading to the fastest convergence rate of consensus. 
Given the wide range of applications of continuous-time consensus \cite{survey}, advancing the understanding of this problem could enhance their effectiveness and broaden their scope. 
However, this open optimization problem is highly challenging and very little is known about it in the existing literature. 
Even with powerful computational resources capable of executing such a combinatorial search, the graphical properties of these optimal directed graphs remain elusive, let alone the associated expensive computational complexity.
To the best of our knowledge, there is no existing literature that directly addresses this problem. The most closely related works are \cite{ming} and \cite{bayen}, which aim to identify optimal directed graphs under specific constraints for other definitions of algebraic connectivity. The special case of the problem -- identifying optimal undirected graphs with maximal algebraic connectivity, given a fixed number of vertices and (undirected) edges -- has been studied in \cite{ogiwara} and \cite{cdc24laplacian}.

This paper aims to address the aforementioned open problem. For sparse and dense graphs, we identify the class of all optimal directed graphs that maximize algebraic connectivity. For general graphs, we propose a computationally efficient sequence of ``almost regular'' directed graphs that achieve fast consensus, with algebraic connectivity close to the optimal possible value.

\section{Optimal Directed Graphs}

We begin with two special cases, namely very sparse and dense graphs, for which the optimal graphs with maximal algebraic connectivity can be explicitly characterized.

To state the first result, we need the following concept.  
A vertex $i$ in a directed graph $\bbb G$ is called a root of $\bbb G$ if for each other vertex $j$ of $\bbb G$, there is a directed path from $i$ to $j$. We say that $\bbb G$ is rooted at vertex $i$ if $i$ is in fact a root, and that $\bbb G$ is rooted if it possesses at least one root. In other words, a directed graph is rooted if it contains a directed spanning tree. 
An $n$-vertex directed tree is a rooted graph with $n-1$ arcs. It is easy to see that a directed tree has a unique root with an in-degree of 0, while all other vertices have an in-degree of exactly 1.
The smallest possible directed tree is a single isolated~vertex.

\vspace{.05in}

\begin{theorem}\label{th:aG_n-1}
Among all simple directed graphs with $n$ vertices and $n-1$ arcs, the maximal algebraic connectivity is $1$, which is achieved if, and only if, the graph is a directed tree.
\end{theorem}

\vspace{.05in}

The theorem implies that when the number of directed edges $m=n-1$, all directed trees achieve the same maximal algebraic connectivity. 
This differs from undirected trees, where the algebraic connectivity varies with the tree structure and reaches its maximum in an undirected star \cite[Theorem 6.5]{Mohar1991}.
To prove Theorem \ref{th:aG_n-1}, we need the following lemmas.

\vspace{.05in}

\begin{lemma}\label{lm:positive}
For any directed graph $\bbb G$, 
    $a(\bbb G)$ is positive if, and only if, $\bbb G$ is rooted. 
\end{lemma}

\vspace{.05in}

{\bf Proof of Lemma \ref{lm:positive}:}
Let $D_{{\rm out}}$ denote the out-degree matrix of $\bbb G$, which is assumed to have $n$ vertices; this matrix is an $n\times n$ diagonal matrix whose $i$th diagonal entry equals the out-degree of vertex $i$. The out-degree based Laplacian matrix is defined as $L_{{\rm out}}=D_{{\rm out}}-A'$. It is straightforward to verify that the (in-degree) Laplacian matrix of a directed graph $\bbb G$ equals the out-degree Laplacian matrix of its transpose graph\footnote{The transpose of a directed graph $\bbb G$ is a directed graph which results when the directed edges in $\bbb G$ are reversed.} $\bbb G'$. In other words, $L(\bbb G)= L_{{\rm out}}(\bbb G')$. Note that the set of all possible simple directed graphs with $n$ vertices is invariant under the graph transpose operation. 
Lemma 2 in \cite{Wu05} shows that the second smallest real part of all eigenvalues of $L_{{\rm out}}(\bbb G)$ is positive if and only if $\bbb G'$ is rooted, which consequently implies that $a(\bbb G)$ is positive if and only if $\bbb G$ is rooted.
\hfill$\qed$

\vspace{.05in}

\begin{lemma}\label{lem:aG-directed-upper}
For any simple directed graph $\bbb G$ with $n$ vertices and $m$ arcs, $a(\bbb G) \le \frac{m}{n-1}$.
\end{lemma}

\vspace{.05in}

Simple examples show that this upper bound may not be achievable. 
As will be seen in Theorem \ref{thm:starsunion}, it is achievable under certain conditions.

\vspace{.05in}



{\bf Proof of Lemma \ref{lem:aG-directed-upper}:}
Label the $n$ Laplacian eigenvalues of $\bbb G$ as $\lambda_1, \lambda_2,\ldots, \lambda_n$ such that ${\rm Re}(\lambda_1)\le {\rm Re}(\lambda_2)\le \cdots \le {\rm Re}(\lambda_n)$. It is clear that $\lambda_1 = 0$.
Note that the $i$th diagonal entry of the Laplacian matrix equals the in-degree of vertex $i$, denoted by $d_i$. 
Since the sum of all eigenvalues of a matrix equals its trace, 
$\sum_{i=1}^n \lambda_i = {\rm tr}(L(\bbb G)) = \sum_{i=1}^n d_i = m$. Thus,
$ m = \sum_{i=2}^n {\rm Re}(\lambda_i) \ge (n-1){\rm Re}(\lambda_2) = (n-1)a(\bbb G)$, 
which implies $a(\bbb G) \le \frac{m}{n-1}$.
\hfill $\qed$

\vspace{.05in}

A directed graph is acyclic if it contains no directed cycles. Thus, by definition, a directed acyclic graph cannot contain a self-arc. Any directed tree is acyclic. The transpose of an acyclic graph remains acyclic.


\begin{lemma}\label{lm:acyclic}
For any acyclic simple directed graph, its Laplacian spectrum consists of its in-degrees.
\end{lemma}

\vspace{.05in}

{\bf Proof of Lemma \ref{lm:acyclic}:}
The adjacency matrix $A$ of a directed graph $\bbb G$, as defined in the introduction, is based on in-degrees. The out-degree based adjacency matrix is the transpose of the in-degree based adjacency matrix; that is, its $ij$th entry equals 1 if $(i,j)$ is an arc in the graph, and equals 0 otherwise, as referenced in \cite[page 151]{Harary69}. 
For any permutation matrix $P$, $P'AP$ represents an adjacency matrix of the same graph, but with its vertices relabeled; the same property applies to out-degree based adjacency matrices.
Since $\bbb G$ is acyclic, from \cite[Theorem 16.3]{Harary69}, there exists a permutation matrix $P$ with which $P'A'P$ is upper triangular. Then, $P'AP$ is lower triangular, which implies that $P'LP$ is also lower triangular. Thus, the spectrum of $P'LP$ consists of its diagonal entries. Since $P'LP$ and $L$ share the same spectrum and diagonal entries, the spectrum of $L$ consists of its diagonal entries, which are the in-degrees of $\bbb G$. 
\hfill $\qed$

\vspace{.05in}

{\bf Proof of Theorem \ref{th:aG_n-1}:}
From Lemma \ref{lem:aG-directed-upper}, for any simple directed graph $\bbb G$ with $n$ vertices and $n-1$ arcs, $a(\bbb G) \le 1$. Thus, to prove the theorem, it is sufficient to show that $a(\bbb G) = 1$ if and only if the graph is a directed tree. 
First, suppose that $a(\bbb G) = 1$. Then, from Lemma \ref{lm:positive}, $\bbb G$ is rooted. Since $\bbb G$ has $n$ vertices and $n-1$ arcs, it must be a directed tree. Conversely, suppose that $\bbb G$ is a directed tree, in which a single vertex has in-degree 0 and $n-1$ vertices have in-degree 1. 
From Lemma \ref{lm:acyclic}, its Laplacian spectrum consists of one eigenvalue at 0 and $n-1$ eigenvalues at 1, which implies that $a(\bbb G) = 1$. 
\hfill $\qed$

\vspace{.05in}

We next consider dense graphs. The following concepts will be needed. 

The union of two directed graphs, $\bbb G_1$ and $\bbb G_2$, with the same vertex set, denoted by $\bbb G_1\cup\bbb G_2$, is the directed graph with the same vertex set and its directed edge set being the union of the directed edge sets of $\bbb G_1$ and $\bbb G_2$. Similarly, the intersection of two directed graphs, $\bbb G_1$ and $\bbb G_2$, with the same vertex set, denoted by $\bbb G_1\cap\bbb G_2$, is the directed graph with the same vertex set and its directed edge set being the intersection of the directed edge sets of $\bbb G_1$ and $\bbb G_2$.
Graph union is an associative binary operation, and thus the definition extends unambiguously to any finite sequence of directed graphs.
The complement of a simple directed graph $\bbb G$, denoted by $\overline{\bbb G}$, is the simple directed graph with the same vertex set such that $\bbb G \cup \overline{\bbb G}$ equals the complete graph and $\bbb G \cap \overline{\bbb G}$ equals the empty graph. 
It is easy to see that if vertex $i$ has in-degree $d_i$ in $\bbb G$, then it has in-degree $n-1-d_i$ in $\overline{\bbb G}$. Moreover, the total number of arcs in $\bbb G$ and $\overline{\bbb G}$ is $n(n-1)$.



The disjoint union of two directed graphs is a larger directed graph whose vertex set is the disjoint union of their vertex sets, and whose arc set is the disjoint union of their arc sets. Disjoint union is an associative binary operation, and thus the definition extends unambiguously to any finite sequence of directed graphs. Any disjoint union of two or more graphs is necessarily disconnected.
A directed forest is a disjoint union of directed tree(s). A directed forest composed of $k$ directed trees thus has $k$ vertices with an in-degree of 0, while all other vertices have an in-degree of exactly 1.
It is easy to see that the number of directed trees in a directed forest is equal to the difference between the number of vertices and the number of arcs.

\vspace{.05in}

\begin{theorem}\label{th:aG-dense}
Among all simple directed graphs with $n$ vertices and $m$ arcs, with $(n-1)^2 \le m < n(n-1)$, the maximal algebraic connectivity is $n-1$, which is achieved if, and only if, the complement of the graph is a directed forest consisting of $m - n(n-2)$ directed trees. 
\end{theorem}

\vspace{.05in}

\begin{figure}[!ht]
\centering
\includegraphics[width=2.2in]{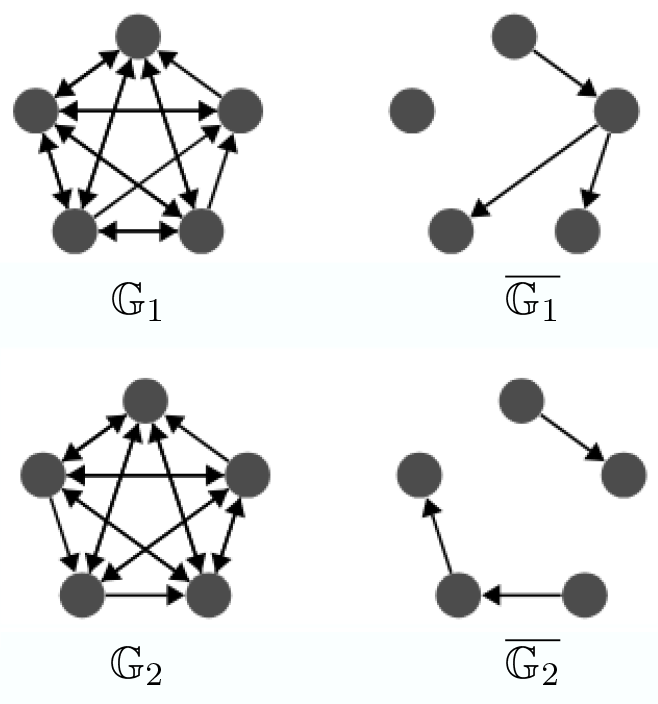} 
\caption{Two examples for Theorem \ref{th:aG-dense} with $n=5$ and $m=17$}
\label{fig:thmcomplement}
\end{figure}

The simple examples in Figure \ref{fig:thmcomplement} serve to illustrate the theorem. For the sake of simplicity in drawing, we use a bidirectional edge to represent two arcs with opposite directions throughout this paper; so each bidirectional edge counts as two arcs. Both $\bbb G_1$ and $\bbb G_2$ are optimal graphs with maximal algebraic connectivity for $n=5$ and $m=17$. The complement of $\bbb G_1$, $\overline{\bbb G_1}$, is a directed forest consisting of $m - n(n-2)=2$ directed trees (one of which is a single isolated vertex). The same observation applies to the complement of $\bbb G_2$, $\overline{\bbb G_2}$, which is consistent with the theorem statement. From the theorem and the examples, there may be multiple optimal graphs with the same maximal algebraic connectivity under the theorem condition.

Theorem \ref{th:aG-dense} does not account for the complete directed graph, which has $n(n-1)$ arcs and an algebraic connectivity of $n$. 
To prove the theorem, we need the following results.

\vspace{0.05in}

\begin{lemma}\label{lm:components}
(Proposition 1.2.11 in \cite{West00})
Any simple undirected graph with $n$ vertices and $m$ edges has at least $\max\{n-m,1\}$ connected components.
\end{lemma}

\vspace{.05in}

We say that a subgraph of a directed graph $\bbb G$ spans $\bbb G$ if it contains all vertices of $\bbb G$.

\vspace{0.05in}


\begin{lemma}\label{lm:multiplicity0}
For any simple directed graph $\bbb G$, the algebraic multiplicity of eigenvalue 0 in its Laplacian spectrum equals the minimum number of directed trees in any directed forest that spans $\bbb G$.
\end{lemma}


\vspace{0.05in}

\begin{figure}[!ht]
\centering
\includegraphics[width=3in]{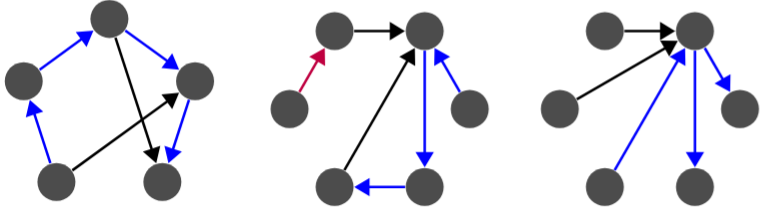} 
\caption{Three examples of spanning directed forests}
\label{fig:span}
\end{figure}

The above lemma generalizes Lemma \ref{lm:positive}. 
Three examples are provided in Figure \ref{fig:span} to illustrate this generalization. 
The Laplacian matrix of the left graph in Figure~\ref{fig:span} has one eigenvalue at 0, and there is one directed tree subgraph that can span the graph, one possible instance of which is depicted in blue. The middle graph in Figure \ref{fig:span} has two eigenvalues at 0 in its Laplacian spectrum, and the minimum number of directed trees in any directed forest that spans the graph is also 2. One such spanning directed forest is depicted in purple and blue, representing its two directed trees, respectively. The right graph in Figure \ref{fig:span} has three eigenvalues equal to 0, and the minimum number of directed trees in any directed forest that spans the graph is also 3. One such spanning forest consists of a directed tree with three vertices depicted in blue, along with two isolated vertices, each considered a trivial tree.

\begin{figure}[!ht]
\centering
\includegraphics[width=3in]{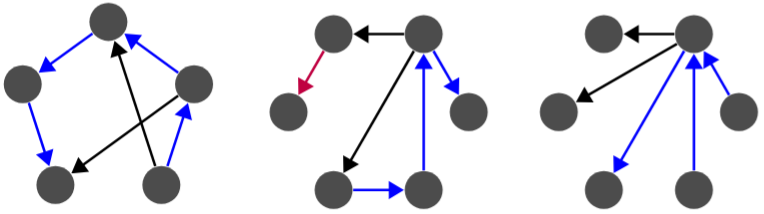} 
\caption{Three examples of spanning in-forests}
\label{fig:intreespan}
\end{figure}

Lemma \ref{lm:multiplicity0} is a direct consequence of Corollary 1 in \cite{Chebotarev02}, which proves that the algebraic multiplicity of the eigenvalue 0 in the out-degree Laplacian spectrum of a simple directed graph equals the minimum number of ``in-trees'' in any ``in-forest'' that spans the graph. The transpose of a directed graph is a directed graph with the same vertex set, but with all arcs reversed in direction compared to the corresponding arcs in the original graph. A directed graph is called an in-tree or in-forest if its transpose is a directed tree or a directed forest, respectively.
Three examples illustrating the corollary are presented in Figure \ref{fig:intreespan}, which are respectively the transposes of the three graphs in Figure \ref{fig:span}.

Return to the statement of Lemma \ref{lm:multiplicity0}. Let $\bbb G'$ be the transpose of graph $\bbb G$. 
It is straightforward to verify that the out-degree based Laplacian matrix (see the definition in the second paragraph of the introduction) of $\bbb G'$ is equal to the (in-degree based) Laplacian matrix of $\bbb G$.
From Corollary 1 in \cite{Chebotarev02}, the algebraic multiplicity of eigenvalue 0 in the out-degree Laplacian spectrum of $\bbb G'$ equals the minimum number of in-trees in any in-forest that spans $\bbb G'$. Therefore, from the relationship between in-trees and directed trees, as well as between in-forests and directed forests, the algebraic multiplicity of the eigenvalue 0 in the (in-degree) Laplacian spectrum of $\bbb G$ equals the minimum number of directed trees in any directed forest spanning $\bbb G$.

We here provide a more direct proof of Lemma \ref{lm:multiplicity0} using a generalized version of the well-known Kirchhoff's matrix tree theorem. To state the generalized matrix tree theorem, we need the following concept.

A weighted simple directed graph is a simple directed graph in which each arc is assigned a nonzero real number as its weight. An unweighted simple directed graph can be viewed as a special case of a weighted graph where all weights are equal to 1. 
The Laplacian matrix of a weighted simple directed graph is defined as the difference between the weighted in-degree matrix and the adjacency matrix, where the weighted in-degree matrix is a diagonal matrix whose $i$th diagonal entry equals the sum of all in-degree arc weights of vertex $i$, and the adjacency matrix is defined such that its $ij$th entry equals the weight of arc $(j,i)$, or 0 if no such arc exists. This definition of the Laplacian clearly simplifies to the one given in the introduction when the graph is unweighted. 
The definition implies that any Laplacian matrix has all row sums equal to 0. It is easy to see that any real matrix with all row sums equal to 0 can be interpreted as the Laplacian matrix of a uniquely determined weighted simple directed graph.

For any simple weighted graph $\bbb G$, let $w(\bbb G)$ denote the product of the weights of all arcs in $\bbb G$. In the special case when $\bbb G$ has no arcs, $w(\bbb G)$ is defined to be 1.

\vspace{.05in}

\begin{lemma}\label{th:matrix-tree-weighted}
{\rm (second Theorem on page 379 of \cite{Chaiken78})}
Let $M$ be any $n \times n$ real matrix with all row sums equal to zero, and let $\bbb G$ denote the weighted simple directed graph uniquely determined by $M$.
For any $1\le k\le n$ distinct indices $i_1, \ldots, i_k\in\{1,\ldots,n\}$, the $(n-k) \times (n-k)$ principal minor of $M$ obtained by removing its rows and columns indexed by $i_1, \ldots, i_k$ is equal to
$\sum_{\bbb F \in \mathcal{F}} w(\bbb F)$,
where $\mathcal{F}$ is the set of spanning directed forests\footnote{A directed forest is termed an arborescence in \cite[page 379]{Chaiken78}. It is implicitly assumed in the proof on \cite[page 380]{Chaiken78} that each directed forest is spanning.} of $\bbb G$ composed of $k$ directed trees rooted at vertices $i_1, \ldots, i_k$.
\end{lemma}

\vspace{.05in}

In the special case when $k=n$, there is only one spanning directed forest composed of $k$ directed trees rooted at $k$ distinct indices; this forest is the spanning subgraph of $\bbb G$ without any arcs, consisting of $n$ isolated vertices. Then, $\sum_{\bbb F \in \mathcal{F}} w(\bbb F)=1$, which is consistent with the convention that a minor of order zero is defined as 1.

In another special case when all arc weights of $\bbb G$ are equal to 1, each $w(\bbb F)$, $\bbb F\in\scr F$ equals 1, and thus $\sum_{\bbb F \in \mathcal{F}} w(\bbb F)$ simplifies to the number of spanning directed forests in $\scr F$. This leads to the following corollary.

\vspace{.05in}

\begin{corollary}\label{th:matrix-tree-unweighted}
Let $L$ be the Laplacian matrix of a simple directed graph $\bbb G$ with $n$ vertices. 
For any $1\le k \le n$ distinct indices $i_1, \ldots, i_k\in\{1,\ldots,n\}$, the $(n-k) \times (n-k)$ principal minor of $L$ obtained by removing its rows and columns indexed by $i_1, \ldots, i_k$ is equal to the number of spanning directed forests of $\bbb G$ composed of $k$ directed trees rooted at vertices $i_1, \ldots, i_k$.
\end{corollary}

\vspace{.05in}

We can now provide a simple proof of Lemma \ref{lm:multiplicity0}.

\vspace{.05in}

{\bf Proof of Lemma \ref{lm:multiplicity0}:} 
Write the characteristic polynomial of $L$ as $p(\lambda) = \det(\lambda I - L)$. It is well known that $p(\lambda) = \sum_{k=0}^n (-1)^{n-k} a_k \lambda^k$, where $a_k$ is the sum of all principal minors of $L$ of order $n-k$. 
In particular, $a_0 = \det (L)=0$ because $L$ has an eigenvalue at 0.
From Corollary \ref{th:matrix-tree-unweighted}, $a_k$ equals the total number of spanning directed forests of $\bbb G$ consisting of exactly $k$ directed trees. 
Let $q$ denote the minimum number of directed trees in any spanning directed forest of $\bbb G$. Then, $a_k=0$ for all $k\in\{1,\ldots,q-1\}$, and $a_k>0$ for all $k\in\{q,\ldots,n\}$. 
From the preceding discussion, $p(\lambda) = \sum_{k=q}^n (-1)^{n-k} a_k \lambda^k = \lambda^q \sum_{k=q}^n (-1)^{n-k} a_k \lambda^{k-q}$, which implies that the algebraic multiplicity of eigenvalue 0 is equal to $q$.
\hfill $\qed$

\vspace{.05in}


We will also need the lemma below.

\vspace{.05in}

\begin{lemma}\label{lm:complement}
If the Laplacian spectrum of an $n$-vertex simple directed graph $\bbb G$ is $\{0,\lambda_2, \ldots, \lambda_n\}$ with $0\le {\rm Re}(\lambda_2)\le \cdots \le {\rm Re}(\lambda_n)$, then the Laplacian spectrum of its complement $\overline{\bbb G}$ is $\{0, n-\lambda_n, \ldots, n-\lambda_2\}$ and $0\le {\rm Re}(n-\lambda_n)\le \cdots \le {\rm Re}(n-\lambda_2)$.
\end{lemma}

\vspace{.05in}

The following proof of the lemma employs the same technique as that used in the proof of Theorem 2 in \cite{Agaev05}, which was developed for a variant of Laplacian matrices.
For any square matrix $M$, we denote its characteristic polynomial as $p_M(\lambda) = \det(\lambda I - M)$ in the sequel.

\vspace{.05in}

{\bf Proof of Lemma \ref{lm:complement}:}
Let $L$ and $\overline{L}$ be the Laplacian matrices of $\bbb G$ and $\overline{\bbb G}$, respectively. It is straightforward to verify that $L+\overline{L} = nI-J$, where $I$ is the identity matrix and $J$ is the $n \times n$ matrix with all entries equal to 1.
Let $Q = L+J = nI-\overline{L}$. 
We first show that 
\begin{equation}\label{eq:complement-pf}
\lambda p_Q(\lambda) = (\lambda-n) p_L(\lambda).
\end{equation}
Note that $p_L(\lambda)=0$ when $\lambda=0$, as $L$ has an eigenvalue at 0. Thus, \eqref{eq:complement-pf} holds when $\lambda=0$.

To prove \eqref{eq:complement-pf} for $\lambda\neq 0$, let $c_i$, $i\in\{1,\ldots,n\}$ denote the $i$th column of matrix $\lambda I - L$. Since $Q = L+J$, it follows that the $i$th column of matrix $\lambda I-Q$ is $c_i-\1$. 
Since the determinant of a matrix is multilinear and adding one column to another does not alter its value, 
$p_Q(\lambda) = \det\; [c_1-\1, c_2-\1, \cdots, c_n-\1] = \det\; [c_1, c_2-\1, \cdots, c_n-\1] - \det\; [\1, c_2-\1, \cdots, c_n-\1] = \det\; [c_1, c_2-\1, \cdots, c_n-\1] - \det\; [\1, c_2, \cdots, c_n]$. Repeating this process sequentially for the columns from 2 to $n$ leads to
\begin{align*}
    p_Q(\lambda) = p_L(\lambda) 
    - \sum_{i=1}^n \det \big[c_1, \cdots, c_{i-1},\1, c_{i+1},\cdots,c_n\big].
\end{align*}
Note that $\sum_{j=1}^n c_j =\lambda\1$, as each row sum of $\lambda I - L$ is equal to $\lambda$. Then, for any $i\in\{1,\ldots,n\}$,
\begin{align*}
    &\det \big[c_1, \cdots, c_{i-1},\1, c_{i+1},\cdots,c_n\big] \\
    =\;& \det \big[c_1, \cdots, c_{i-1},\textstyle\frac{1}{\lambda}\sum_{j=1}^n c_j, c_{i+1},\cdots,c_n\big] \\
    =\;& \textstyle\frac{1}{\lambda}\det \big[c_1, \cdots, c_{i-1},\sum_{j=1}^n c_j, c_{i+1},\cdots,c_n\big] \\
    =\;& \textstyle\frac{1}{\lambda}\det \big[c_1, \cdots, c_n\big] = \frac{1}{\lambda} p_L(\lambda).
\end{align*}
Substituting this equality into the preceding expression for $p_Q(\lambda)$ yields $p_Q(\lambda) = p_L(\lambda) - \frac{n}{\lambda}p_L(\lambda)$, which proves \eqref{eq:complement-pf}.

To proceed, recall that $\overline L = nI-Q$. Then, 
\begin{align*}
p_{\overline{L}}(\lambda) &= \det(\lambda I - \overline{L}) = \det(\lambda I - nI + Q) \\
&= (-1)^n \det((n-\lambda)I - Q) = (-1)^n p_Q(n-\lambda). 
\end{align*}
From this and \eqref{eq:complement-pf}, with $\lambda$ substituted by $n-\lambda$,
\begin{align}\label{eq:eigrelation}
    (n-\lambda)p_{\overline{L}}(\lambda) = (-1)^{n+1} \lambda p_L(n-\lambda).
\end{align}
Both sides of \eqref{eq:eigrelation} are polynomials in $\lambda$ of degree $n+1$. It is easy to see that 0 and $n$ are roots of both sides, as 0 is an eigenvalue of both $L$ and $\overline{L}$. Then, the set of nonzero roots of $p_{\overline{L}}(\lambda)$ coincides with the set of roots of $p_L(n-\lambda)$, excluding the root at $n$.
Therefore, if the Laplacian spectrum of $\bbb G$ is $\{0,\lambda_2, \ldots, \lambda_n\}$, then the Laplacian spectrum of $\overline{\bbb G}$ is $\{0, n-\lambda_n, \ldots, n-\lambda_2\}$. Recall that the smallest real part of all Laplacian eigenvalues is always 0. With these facts, it is easy to see that if $0\le {\rm Re}(\lambda_2)\le \cdots \le {\rm Re}(\lambda_n)$, then $0\le {\rm Re}(n-\lambda_n)\le \cdots \le {\rm Re}(n-\lambda_2)$.
\hfill $\qed$


\vspace{.05in}

We are now in a position to prove Theorem \ref{th:aG-dense}.

\vspace{.05in}

{\bf Proof of Theorem \ref{th:aG-dense}:}
Let the Laplacian spectrum of an $n$-vertex simple directed graph $\bbb G$ be $\{0,\lambda_2, \ldots, \lambda_n\}$ with $0\le {\rm Re}(\lambda_2)\le \cdots \le {\rm Re}(\lambda_n)$ and the Laplacian spectrum of its complement $\overline{\bbb G}$ be $\{0,\mu_2, \ldots, \mu_n\}$ with $0\le {\rm Re}(\mu_2)\le \cdots \le {\rm Re}(\mu_n)$. From Lemma \ref{lm:complement}, $\mu_i = n-\lambda_{n+2-i}$ for each $i \in \{2,\ldots,n\}$. 

Let $\bbb H$ be the underlying simple undirected graph of $\overline{\bbb G}$, obtained by replacing all directed edges in $\overline{\bbb G}$ with undirected ones. It is easy to see that $\bbb H$ has $n(n-1)-m$ (undirected) edges. From Lemma \ref{lm:components}, $\bbb H$ has at least $n-(n(n-1)-m) = m-n(n-2)$ connected components. Thus, $\overline{\bbb G}$ has at least $m-n(n-2)$ weakly connected components. 

To span a weakly connected component, at least one directed tree is required. Thus, the minimum number of directed trees required to form a directed forest that spans $\overline{\bbb G}$ is at least $m-n(n-2)$. From Lemma \ref{lm:multiplicity0}, the algebraic multiplicity of the eigenvalue 0 in the Laplacian spectrum of $\overline{\bbb G}$ is at least $m-n(n-2)$. 
Since the smallest real part of all Laplacian eigenvalues is always 0, $\mu_i=0$ for each $i\in\{2,\ldots,m-n(n-2)\}$.
As the trace of a square matrix equals the sum of its eigenvalues, ${\rm tr}(\overline{L}) = \sum_{i=m-n(n-2)+1}^{n} \mu_i $, where $\overline{L}$ denotes the Laplacian matrix of $\overline{\bbb G}$.
Meanwhile, as the $i$th diagonal entry of $\overline{L}$ is $n-1-d_i$, where $d_i$ denotes the in-degree of vertex $i$ in $\bbb G$, ${\rm tr}(\overline{L}) = \sum_{i=1}^n (n-1-d_i) = n(n-1)-m$. It follows that
\begin{align}
    n(n-1)-m &= \textstyle\sum_{i=m-n(n-2)+1}^{n} {\rm Re}(\mu_i) \nonumber\\
    &\le (n(n-1)-m)  {\rm Re}(\mu_n),\label{eq:thm4-2}
\end{align}
which implies that ${\rm Re}(\mu_n) \ge 1$, and thus 
$a(\bbb G) = {\rm Re}(\lambda_2) = n-{\rm Re}(\mu_n)\le n-1$.

We first prove the sufficiency of the theorem. Suppose that $a(\bbb G) = n-1$. It follows that ${\rm Re}(\mu_n)=1$, and thus the equality in \eqref{eq:thm4-2} must hold, which implies that ${\rm Re}(\mu_i) = {\rm Re}(\mu_n)=1$ for all $i \in \{m-n(n-2)+1,\ldots, n\}$. 
Then, the algebraic multiplicity of eigenvalue 0 of $\overline{L}$ is exactly $m-n(n-2)$. From Lemma \ref{lm:multiplicity0}, the minimum number of directed trees required to form a directed forest that spans $\overline{\bbb G}$ is $m-n(n-2)$. Recall that $\overline{\bbb G}$ has at least $m-n(n-2)$ weakly connected components. It follows that $\overline{\bbb G}$ has exactly $m-n(n-2)$ weakly connected components, and thus each of them can be spanned by a directed tree. 
It is easy to see that a weakly connected graph can be spanned by a directed tree if and only if the graph is rooted. Then, each of the $m-n(n-2)$ weakly connected components is rooted. Without loss of generality, label these weakly connected components from 1 to $m-n(n-2)$. Let $n_i$ denote the number of vertices in component $i$. It is clear that $\sum_{i=1}^{m-n(n-2)}n_i=n$. Since each component is rooted, component~$i$ has at least $n_i-1$ arcs. Then, the total number of arcs in $\overline{\bbb G}$ is at least $\sum_{i=1}^{m-n(n-2)} (n_i-1) = n - (m-n(n-2)) = n(n-1)-m$. As $\bbb G$ has $m$ arcs, $\overline{\bbb G}$ has exactly $n(n-1)-m$ arcs. It follows that each component $i$ must have exactly $n_i-1$ arcs, which implies that each component is a directed tree. Therefore, $\overline{\bbb G}$ is a directed forest consisting of $m - n(n-2)$ directed trees.



We next prove the necessity of the theorem. Suppose that $\overline{\bbb G}$ is a directed forest consisting of $m - n(n-2)$ directed trees. It is easy to see that the in-degree sequence of $\overline{\bbb G}$ consists of $m - n(n-2)$ zeros and $n(n-1)-m$ ones. Since $\overline{\bbb G}$ is acyclic, from Lemma \ref{lm:acyclic}, the Laplacian spectrum of $\overline{\bbb G}$ coincides with its in-degree sequence. It follows that $\mu_n = 1$, and thus $a(\bbb G) = {\rm Re}(\lambda_2) = {\rm Re}(n-\mu_n) = n-1$.
\hfill $\qed$

\vspace{.05in}

For general cases when the number of vertices $n$ and the number of arcs $m$ satisfy $n-1 < m < (n-1)^2$, theoretically identifying the optimal graphs with maximal algebraic connectivity remains challenging and has thus far eluded us, except in the following special cases.

\vspace{.05in}

\begin{theorem}\label{thm:starsunion}
Among all simple directed graphs with $n$ vertices and $m$ arcs, with $m = l(n-1)$ and $l\in\{1,\ldots,n\}$, the maximal algebraic connectivity is $l$, which is achieved if the graph is the union of $l$ simple $n$-vertex directed stars, each rooted at a distinct vertex. 
\end{theorem}

\vspace{.05in}

\begin{figure}[!ht]
\centering
\includegraphics[width=2in]{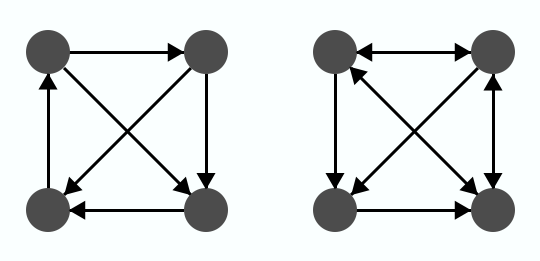} 
\caption{Two $4$-vertex optimal graphs respectively with $6$ and $9$ arcs}
\label{fig:notstarsunion}
\end{figure}

The above theorem is a direct consequence of Lemma \ref{lem:aG-directed-upper} and Theorem \ref{th:p1-construction}, which will be discussed in the next section after the statement of Theorem \ref{th:p1-construction}. 
It is worth emphasizing that the graphical condition identified in Theorem \ref{thm:starsunion} is not necessary. A simple example arises when $m=n-1$, where all other optimal graphs have been identified in Theorem \ref{th:aG_n-1}. Figure \ref{fig:notstarsunion} illustrates two optimal graphs, in which the number of arcs $m$ is a multiple of $n-1$, yet they are not a union of directed stars.

It is also worth noting that the classes of optimal graphs identified in Theorems \ref{th:aG_n-1} to \ref{thm:starsunion} are all ``almost regular'' directed~graphs.

\section{Almost Regular Directed Graphs}

A directed graph is called almost regular if the difference between its largest and smallest in-degrees is at most 1.
It is easy to see that any directed forest, including a directed tree, is almost regular, and so are the optimal graphs identified in Theorem \ref{th:aG_n-1}. It is also evident that any union of directed stars of the same size, rooted at distinct vertices, as identified in Theorem \ref{thm:starsunion}, is almost regular. 

Theorem \ref{th:aG-dense} states that the complement of any optimal graph is a directed forest consisting of $m - n(n-2)$ directed trees, which implies that the complement graph has $m - n(n-2)$ vertices with in-degree 0, while all other vertices have in-degree 1. Then, in the optimal graph, those $m - n(n-2)$ vertices have in-degree $n-1$, and all the other vertices have in-degree $n-2$. Therefore, all the optimal graphs identified in Theorem \ref{th:aG-dense} are almost regular.

\begin{figure}[!ht]
\centering
\includegraphics[width=2.2in]{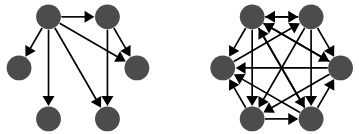} 
\caption{Two optimal graphs with maximal algebraic connectivity that are not almost regular}
\label{fig:notregular}
\end{figure}

The following examples demonstrate that there exist optimal graphs with maximal algebraic connectivity that are not almost regular. Specifically, the left graph in Figure \ref{fig:notregular} is an optimal graph with maximal algebraic connectivity, given 6 vertices and 7 arcs, while the right graph has maximal algebraic connectivity with 6 vertices and 22 arcs; neither graph is almost regular.

In the following, we present an algorithm that inductively constructs a sequence of graphs with $n$ vertices, each of which is a simple directed graph with (almost) maximal algebraic connectivity for every possible number of arcs.

The algorithm will make use of the following notation. We use $\lceil \cdot \rceil$ to denote the ceiling function, which maps a real number $x$ to the smallest integer greater than or equal to $x$, and use $a\bmod b$ to denote the modulo operation of two integers $a$ and $b$, which returns the remainder after dividing $a$ by $b$. 

\vspace{.05in}

{\em Algorithm 1:} Given $n>1$ vertices, label them, without loss of generality, from $1$ to $n$. Let $\bbb G(n,m)$ denote the $n$-vertex simple directed graph to be constructed with $m$ arcs. Start with the $m=0$ case and set $\bbb G(n,0)$ as the empty graph. For each integer $1\le m\le n(n-1)$, construct $\bbb G(n,m)$ by adding an arc from vertex $\lceil \frac{m}{n-1} \rceil$ to vertex $n-((m-1) \bmod n)$ on $\bbb G(n,m-1)$.
\hfill$\Box$

\vspace{.05in}

Before proceeding, we show that the last sentence in the description of Algorithm 1 is always feasible. To this end, it is sufficient to prove the following two lemmas. Lemma~\ref{lem:validarc} ensures that the arc to be added, $(\lceil \frac{m}{n-1} \rceil, n-((m-1) \bmod n))$, is always a valid arc, and Lemma \ref{lem:algorithmwelldefined} guarantees that adding this arc is always possible. 

\vspace{.05in}

\begin{lemma}\label{lem:validarc}
    For a fixed $n>1$ and any $1\le m\le n(n-1)$, $\lceil \frac{m}{n-1} \rceil$ and $n-((m-1) \bmod n)$ are two distinct integers in $\{1,\ldots,n\}$. 
\end{lemma}

\vspace{.05in}

{\bf Proof of Lemma \ref{lem:validarc}:}
It is clear that $\lceil \frac{m}{n-1} \rceil\ge 1$. Since $m \le n(n-1)$, $\lceil \frac{m}{n-1} \rceil \le \lceil \frac{n(n-1)}{n-1} \rceil = n$. Meanwhile, it is easy to see that $1 \le n-((m-1) \bmod n) \le n$. Thus, both $\lceil \frac{m}{n-1} \rceil$ and $n-((m-1) \bmod n)$ are integers in $\{1,\ldots,n\}$. 
To prove that they are distinct, suppose to the contrary that $\lceil \frac{m}{n-1} \rceil = n-((m-1) \bmod n)$.
Let $q = \lceil \frac{m}{n-1} \rceil$. Then, $(m-1) \bmod n = n-q = (n-q) \bmod n = (-q) \bmod n$, implying $m \equiv -q+1\! \pmod n$.
Meanwhile, $(q-1)(n-1)<m\le q(n-1)$, which implies that $m = q(n-1)-r$ for some integer $r \in \{0,\ldots,n-2\}$. Taking both sides of this equation modulo $n$ leads to 
$m \equiv -q-r\! \pmod n$.
It follows from the two congruence relations above that $r+1 \equiv 0\! \pmod n$.
But this is in contradiction with $r \in \{0,\ldots,n-2\}$. Therefore, $\lceil \frac{m}{n-1} \rceil \neq n-((m-1) \bmod n)$.
\hfill$\qed$


\vspace{.05in}

\begin{lemma}\label{lem:algorithmwelldefined}
    For a fixed $n>1$ and any $1\le m\le n(n-1)$, the arc added during the construction of $\bbb G(n,m)$, $(\lceil \frac{m}{n-1} \rceil, n-((m-1) \bmod n))$, is distinct from all arcs in $\bbb G(n,m-1)$.

\end{lemma}

\vspace{.05in}

{\bf Proof of Lemma \ref{lem:algorithmwelldefined}:}
The lemma will be proved by induction on $m$. 
In the base case when $m=1$, the lemma is clearly true. Now, for the inductive step, suppose that the lemma holds for all $m\in\{1,\ldots,k\}$, where $k$ is a positive integer satisfying $k<n(n-1)$. This implies that all $k$ (distinct) arcs in $\bbb G(n,k)$ are
\eq{\textstyle\big(\big\lceil \frac{i}{n-1} \big\rceil, \; n-((i-1) \bmod n)\big), \;\;\; i \in \{1,\ldots,k\}.\label{eq:allarcs}}
Let $m=k+1$. The corresponding added arc becomes $(\lceil \frac{k+1}{n-1} \rceil, n-(k \bmod n))$. 
To prove that the added arc is distinct from any arc in \eqref{eq:allarcs}, we divide the arcs in \eqref{eq:allarcs} into two categories based on index $i\in \{1,\ldots,k\}$.

First, if $\lceil \frac{i}{n-1} \rceil\neq \lceil \frac{k+1}{n-1} \rceil$, then the added arc clearly differs from $(\lceil \frac{i}{n-1}\rceil, n-((i-1) \bmod n))$. 
Second, if $\lceil \frac{i}{n-1} \rceil = \lceil \frac{k+1}{n-1} \rceil$, then the values of $k+1$ and $i$ differ by at most $n-2$. To see this, set $p=\lceil \frac{k+1}{n-1} \rceil = \lceil \frac{i}{n-1} \rceil$, and thus $p-1<\lceil \frac{k+1}{n-1} \rceil = \lceil \frac{i}{n-1} \rceil$. It follows that both $k+1$ and $i$ lie within the interval $((p-1)(n-1),p(n-1)]$. Since both $k+1$ and $i$ are integers, their difference $k+1-i$ is at most $p(n-1)-(p-1)(n-1)-1=n-2$. With this fact, we further claim that $k+1 \not\equiv i\! \pmod n$. To prove the claim, suppose to the contrary that $k+1 \equiv i\! \pmod n$, which implies that $k+1-i=qn$ for some integer $q$. As $k+1>i$, $q$ is at least 1. But this contradicts the fact that $k+1-i\le n-2$. Thus, the claim is true, which implies that $n-(k \bmod n) \neq n-((i-1) \bmod n)$. Therefore, the added arc is different from $(\lceil \frac{i}{n-1}\rceil, n-((i-1) \bmod n))$. 
The above discussion on the two categories collectively implies that the added arc does not exist in $\bbb G(n,k)$, which completes the inductive step. By induction, the lemma is proved. 
\hfill$\qed$


\vspace{.05in}

\begin{figure}[!ht]
\centering
\includegraphics[width=3.4in]{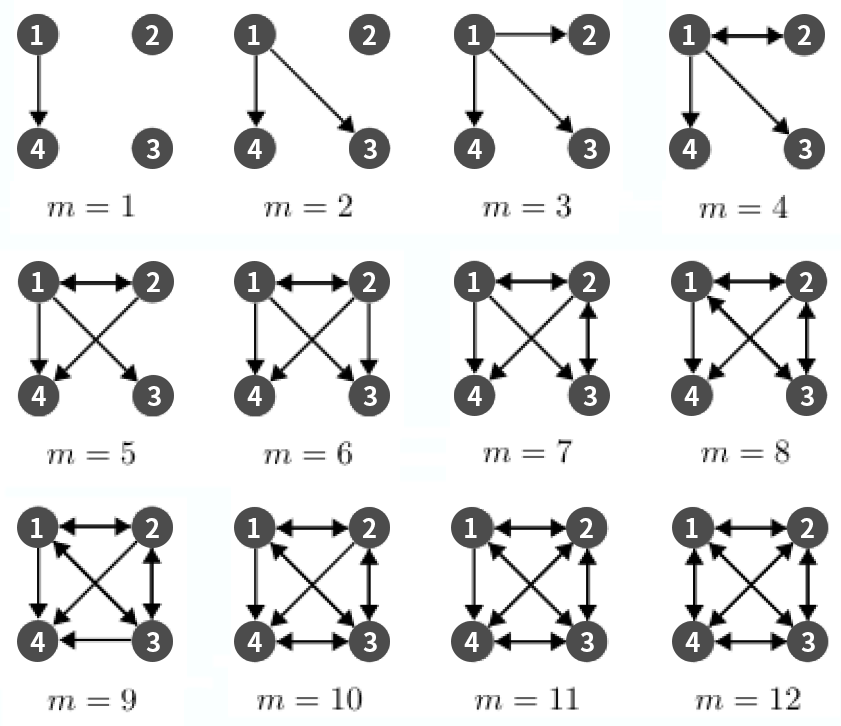} 
\caption{Inductive construction process of Algorithm 1 for $n=4$}
\label{fig:algorithm4nodes}
\end{figure}

Lemmas \ref{lem:validarc} and \ref{lem:algorithmwelldefined} ensure that
Algorithm 1 is properly formulated and operates unambiguously under the given conditions for $n$ and $m$. 
Figure \ref{fig:algorithm4nodes} illustrates the inductive construction process described by the algorithm for the case when $n=4$.
Moreover, from the algorithm description, Lemma \ref{lem:validarc}, and Lemma \ref{lem:algorithmwelldefined}, the arc set of $\bbb G(n,m)$ is composed of $m$ (distinct) arcs of the form 
\eq{\textstyle\big(\big\lceil \frac{i}{n-1} \big\rceil, \; n-((i-1) \bmod n)\big), \;\;\; i \in \{1,\ldots,m\}.\label{eq:allGarcs}}
We will use this fact without special mention in the sequel.

To construct a graph with $n$ vertices and $m$ arcs from scratch using Algorithm 1, the computational complexity is $O(m)$, as identifying the endpoints for each of the $m$ arcs takes $O(1)$ time. 

In the sequel, we set $\kappa = \lfloor \frac{m}{n-1} \rfloor$, where $\lfloor \cdot \rfloor$ denotes the floor function, which maps a real number $x$ to the greatest integer less than or equal to $x$. 
The following result characterizes the graphs constructed by Algorithm 1.

\vspace{.05in}

\begin{proposition}\label{prop:p1-stars}
The graph constructed by Algorithm 1 with $n$ vertices and $m$ arcs, $\bbb G(n,m)$, is the union of $\kappa$ simple $n$-vertex directed stars, rooted at vertices $1, \ldots, \kappa$, respectively, and the $n$-vertex simple directed forest consisting of an $(m-\kappa(n-1)+1)$-vertex directed star, rooted at vertex $\kappa+1$, and $n-(m-\kappa(n-1)+1)$ isolated vertices.
\end{proposition}

\vspace{.05in}

\begin{figure}[!ht]
\centering
\includegraphics[width=2.2in]{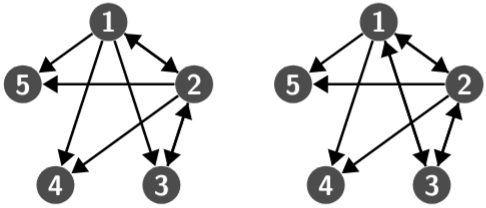} 
\caption{Graphs $\bbb G(5,9)$ and $\bbb G(5,10)$ constructed by Algorithm 1}
\label{fig:prop1-ex}
\end{figure}

Each graph in Figure \ref{fig:algorithm4nodes} validates the 
proposition. We provide some additional explanation for the special case when $m=n-1$, which implies $\kappa=1$. In this case, $\bbb G(n,m)$ in the proposition statement simplifies to the union of one $n$-vertex simple directed star, rooted at vertex $1$, and the $n$-vertex simple directed forest consisting of one $1$-vertex directed star, rooted at vertex $2$, and $n-1$ isolated vertices. Note that a $1$-vertex directed star is simply a single isolated vertex, which implies that the latter graph comprising the union is an empty graph. Therefore, $\bbb G(n,m)$ is exactly an $n$-vertex simple directed star when $m=n-1$; the third graph in Figure \ref{fig:algorithm4nodes} validates this. 
We further provide two 5-vertex examples in Figure \ref{fig:prop1-ex} to illustrate the proposition for general cases. The left graph $\bbb G(5,9)$ is the union of two simple 5-vertex directed stars, respectively rooted at vertices 1 and 2, and a 5-vertex simple directed forest consisting of a 2-vertex directed star rooted at vertex 3 and three isolated vertices indexed as 1, 4, and 5.
Similarly, the right graph $\bbb G(5,10)$ is the union of two simple 5-vertex directed stars, rooted at vertices 1 and 2, respectively, and a 5-vertex directed forest consisting of a 3-vertex directed star rooted at vertex 3 and two isolated vertices labeled 4 and 5.


Algorithm 1 constructs graphs for all possible numbers of arcs. When applying the constructed graphs to consensus, it is understood that $m \ge n-1$, as Proposition \ref{prop:p1-stars} ensures rooted graphs; conversely, if $m<n-1$, the constructed graphs are disconnected.

\vspace{.05in}

{\bf Proof of Proposition \ref{prop:p1-stars}:}
Note that $m$ can be written as $m = \kappa(n-1)+r$, where $0 \le r \le n-2$ is the unique remainder when $m$ is divided by $n-1$.

First, in the case when $r=0$, which implies that $n-1$ divides $m$, all arcs in \eqref{eq:allGarcs} originate from a vertex with an index in the range from 1 to $\lceil \frac{m}{n-1} \rceil = \kappa$. Moreover, each vertex $j\in\{1,\ldots,\kappa\}$ has exactly out-degree $n-1$. 
Thus, the graph is the union of $\kappa$ simple $n$-vertex directed stars, rooted at vertices $1, \ldots, \kappa$, respectively. Note that in this case $m-\kappa(n-1)+1 = 1$. Then, ``the $n$-vertex simple directed forest consisting of an $(m-\kappa(n-1)+1)$-vertex directed star, rooted at vertex $\kappa+1$, and $n-(m-\kappa(n-1)+1)$ isolated vertices'' as described in the proposition statement, reduces to an empty graph. Therefore, the proposition holds in this case.

Next, consider the case when $r \in \{1,\ldots, n-2\}$, which implies that  $\lceil \frac{m}{n-1} \rceil= \kappa+1$. 
Note that in this case each vertex $j\in\{1,\ldots,\kappa\}$ still has exactly out-degree $n-1$, forming $\kappa$ simple $n$-vertex directed stars. 
These $\kappa$ stars consists of $\kappa(n-1)$ arcs. From \eqref{eq:allGarcs}, all remaining $m-\kappa(n-1)$  arcs originate from vertex $\kappa+1$, constituting an $(m-\kappa(n-1)+1)$-vertex directed star. Then, the proposition is clearly true in this case.
\hfill $\qed$

\vspace{.05in}

More can be said about the vertex in-degree sequence of $\bbb G(n,m)$. In the sequel, we set $\nu = \lfloor \frac{m}{n} \rfloor$.



\vspace{.05in}

\begin{proposition}\label{prop:p1-deg}
The graph constructed by Algorithm 1 with $n$ vertices and $m$ arcs, $\bbb G(n,m)$, is almost regular, with $n(\nu+1)-m$ vertices of in-degree $\nu$ and $m-n\nu$ vertices of in-degree $\nu+1$. 
To be more precise, the sequence of its vertex in-degrees is 
\eq{(d_1, \ldots, d_n) = (\underbrace{\nu,\ldots,\nu}_{n(\nu+1)-m}, \; \underbrace{\nu+1,\ldots,\nu+1}_{m-n\nu}\;).\label{eq:indegree}}
\end{proposition}

\vspace{.1in}

{\bf Proof of Proposition \ref{prop:p1-deg}:}
From the expression of the head index of each arc in \eqref{eq:allGarcs}, $n-((i-1) \bmod n)$, 
it takes the value $n+1-i$ for each $i\in\{1,\ldots,n\}$, and then, if $m>n$, repeats with a period of $n$ as $i$ ranges from $n+1$ to $m$. It follows that $d_n \ge \cdots \ge d_1$ and $d_n-d_1\le 1$. Thus, $\bbb G(n,m)$ is almost regular. 
The remaining statement of the proposition directly follows from the following lemma.
\hfill $\qed$


\begin{lemma}\label{lem:degreesequence}
    For any almost regular simple directed graph with $n$ vertices and $m$ arcs, assume, without loss of generality, that its vertex in-degrees satisfy $d_1\le \cdots \le d_n$. Then, its in-degree sequence is \eqref{eq:indegree}.
\end{lemma}

\vspace{.05in}

{\bf Proof of Lemma \ref{lem:degreesequence}:}
Since the graph is almost regular, $d_n-d_1\le 1$. 
Suppose there are $1\le p\le n$ vertices with the minimal in-degree $d_1$. Then, the remaining $q=n-p$ vertices have an in-degree of $d_1+1$. It follows that $m=pd_1+q(d_1+1)=d_1n+q$. As $q$ takes a value in $\{0,1,\ldots,n-1\}$, $d_1$ and $q$ are respectively the unique quotient and remainder when $m$ is divided by $n$. Then, $d_1=\lfloor \frac{m}{n} \rfloor = \nu$ and $q=m-d_1n=m-n\nu$. Therefore, the in-degree sequence is \eqref{eq:indegree}.
\hfill $\qed$


\vspace{.05in}

From Proposition \ref{prop:p1-deg},  the in-degree $d_i$ of each vertex $i$ in $\bbb G(n,m)$ may take an integer value ranging from $0$ to $n-1$, depending on the value of $m$. More observation can be made.

\vspace{.05in}

\begin{lemma}\label{lem:arctail}
For each vertex $i$ in $\bbb G(n,m)$, its $d_i$ incoming arcs originate from $d_i$ vertices whose indices are the $d_i$ smallest elements of $\{1,\ldots,n\}\setminus\{i\}$.
\end{lemma}

\vspace{.05in}

{\bf Proof of Lemma \ref{lem:arctail}:}
From Proposition \ref{prop:p1-stars}, $\bbb G(n,m)$ is composed of $\kappa$ simple $n$-vertex directed stars, rooted at vertices $1, \ldots, \kappa$, respectively, and an $(m-\kappa(n-1)+1)$-vertex directed star rooted at vertex $\kappa+1$. 
The following observations can thus be summarized. First, for each vertex $i \in \{1, \ldots, \kappa\}$, the set of tail indices for all its incoming arcs is either $\{1, \ldots, \kappa\} \setminus \{i\}$ or $\{1, \ldots, \kappa+1\} \setminus \{i\}$. Second, for vertex $\kappa+1$, the set of tail indices for all its incoming arcs is $\{1, \ldots, \kappa\}$. Last, for each vertex $i \in \{\kappa+2, \ldots, n\}$, the set of tail indices for all its incoming arcs is either $\{1, \ldots, \kappa\}$ or $\{1, \ldots, \kappa+1\}$. The lemma therefore follows directly from these observations.
\hfill$\qed$

\vspace{.05in}

From the proof of Lemma \ref{lem:arctail}, it appears that the vertices in $\bbb G(n,m)$ may have three different values of in-degrees when $\kappa\ge 1$, namely $\kappa-1$, $\kappa$, and $\kappa+1$, which seems to conflict with \eqref{eq:indegree} in Proposition \ref{prop:p1-deg} at first glance. Hence, we show here that these two ways of characterizing the in-degree sequence, respectively using $\kappa$ and $\nu$, are consistent. To this end, we take a closer look at the proof of Lemma~\ref{lem:arctail}. The $\kappa$ simple $n$-vertex directed stars contribute either $\kappa-1$ or $\kappa$ to each vertex's in-degree. A vertex with in-degree $\kappa+1$ can exist only if the $(m-\kappa(n-1)+1)$-vertex directed star rooted at vertex $\kappa+1$ contains at least one arc. This condition implies that $n-1$ does not divide $m$, and consequently, $\kappa<n$. As all arcs in the star originate from the vertex $\kappa+1$, from \eqref{eq:allGarcs}, these arcs correspond to the indices $i\in\{1,\ldots,m\}$ such that $\lceil \frac{i}{n-1} \rceil = \kappa+1$. Such indices can be written as $i=\kappa (n-1)+j$, $j\in\{1,\ldots,r\}$, where $r$ is the is the unique remainder when $m$ is divided by $n-1$. Then, the head indices of these arcs are $n-((i-1) \bmod n)=n-((\kappa n-\kappa +j-1) \bmod n)=n-((n-\kappa +j-1) \bmod n)$, $j\in\{1,\ldots,r\}$. 
Note that, since $1\le \kappa\le n-1$, the expression $n-((n-\kappa +j-1) \bmod n)$ takes the value $\kappa+1-j$ if $1\le j\le \kappa$, and evaluates to $n+\kappa+1-j$ if $\kappa+1\le j\le r<n-1$. That is, the head index takes values from $\kappa$ to $1$ as $j$ ranges from $1$ to $\kappa$, and from $n$ to $n+\kappa+1-r>\kappa$ as $j$ ranges from $\kappa+1$ to $r$. From the description of Algorithm 1, the construction order of these arcs follows the ascending order of indices $j\in\{1,\ldots,r\}$. This construction order ensures that vertices with in-degrees $\kappa-1$ and $\kappa+1$ cannot appear simultaneously. The construction process for the $n=4$ case with all possible nonzero $m$, as shown in Figure \ref{fig:algorithm4nodes}, illustrates the above  arguments and conclusion. Therefore, $\bbb G(n,m)$ is almost regular, and by Lemma \ref{lem:degreesequence}, its in-degree sequence is uniquely determined.



Proposition \ref{prop:p1-deg} and Lemma \ref{lem:arctail} together fully characterize the arc set of $\bbb G(n,m)$ and lead to the following consequence.

\vspace{.05in}

\begin{corollary}\label{coro:subgraph}
Let $n_i$ and $m_i$ be integers such that $n_i\ge 2$ and $0\le m_i\le n_i(n_i-1)$, where $i\in\{1,2\}$.
For each~$i$, let $q_i$ and $r_i$ respectively denote the unique quotient and remainder when $m_i$ is divided by $n_i$. If $n_1\le n_2$, $q_1 \le q_2$, and $n_1-r_1 \ge n_2-r_2$, then $\bbb G(n_1,m_1)$ is a subgraph of $\bbb G(n_2,m_2)$.     
\end{corollary}

\vspace{.05in}

{\bf Proof of Corollary \ref{coro:subgraph}:}
Note that $q_i=\lfloor\frac{m_i}{n_i}\!\rfloor$ for each $i\in\{1,2\}$, and $n_1-r_1 \ge n_2-r_2$ implies $n_1(q_1+1)-m_1 \ge n_2(q_2+1)-m_2$. These facts and \eqref{eq:indegree}, along with $q_1 \le q_2$, imply that for each vertex $i\in\{1,\ldots,n_1\}$, its in-degree in $\bbb G(n_1,m_1)$ is no larger than its in-degree in $\bbb G(n_2,m_2)$. It then follows from Lemma \ref{lem:arctail} that each arc in $\bbb G(n_1,m_1)$ is also an arc in $\bbb G(n_2,m_2)$, which completes the proof. 
\hfill$\qed$

\vspace{.05in}

Note that the conditions in Corollary \ref{coro:subgraph} imply $m_1=q_1n_1+r_1\le q_2n_2 + r_1+n_2-n_1\le q_2n_2+r_2 = m_2$. 
In the special case when $n_1=n_2$, from the description of Algorithm 1, $m_1\le m_2$ guarantees that $\bbb G(n_1,m_1)$ is a subgraph of $\bbb G(n_2,m_2)$. In general cases when $n_1<n_2$, simple examples show that $m_1\le m_2$ cannot guarantee this.







The most important property of the graphs constructed by Algorithm 1 is stated in the following theorem.

\vspace{.05in}

\begin{theorem}\label{th:p1-construction}
The algebraic connectivity of the graph constructed by Algorithm 1 with $n$ vertices and $m$ arcs, $\bbb G(n,m)$, equals $\lfloor \frac{m}{n-1} \rfloor$.
\end{theorem}

\vspace{.05in}

Recall that $a(\bbb G) \le \frac{m}{n-1}$ (cf. Lemma \ref{lem:aG-directed-upper}). The above theorem immediately implies that whenever $m$ is a multiple of $n-1$, the graph constructed by Algorithm 1 achieves the maximum algebraic connectivity. From Proposition \ref{prop:p1-stars}, in this case, the constructed graph is the union of multiple directed stars rooted at distinct vertices. This proves Theorem~\ref{thm:starsunion}.
Graphs constructed by Algorithm 1 may also achieve the exact maximum algebraic connectivity when $m$ is not a multiple of $n-1$. For example, $\bbb G(5,9)$, shown as the left graph in Figure \ref{fig:prop1-ex}, has an algebraic connectivity of 2, which is the maximum achievable among all simple directed graphs with 5 vertices and 9 arcs, as identified by an exhaustive simulation search.
Some other examples of small-sized graphs include cases when $n=4$ with $m \in\{3,4\}$, $n=5$ with $m \in\{4,5,6\}$, and $n=6$ with $m \in\{5,6,7\}$.

\begin{figure}[!ht]
\centering
\includegraphics[width=2.2in]{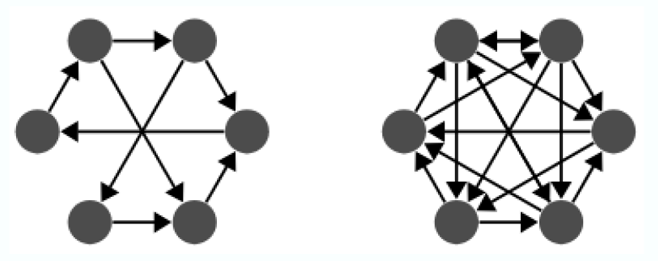} 
\caption{Two optimal graphs with larger algebraic connectivity than the corresponding graphs constructed by Algorithm 1}
\label{fig:almostmaximal}
\end{figure}

Theorem \ref{th:p1-construction} implies that the algebraic connectivity of the graph constructed by Algorithm 1 is generally ``close to'' the maximum possible value, with a gap of no more than 1. For instance, the left graph in Figure \ref{fig:almostmaximal} is an optimal graph with 6 vertices and 8 arcs, exhibiting an algebraic connectivity of $1.123$. In comparison, the graph constructed by Algorithm~1, $\bbb G(6,8)$, has an algebraic connectivity of 1. Similarly, the right graph in Figure \ref{fig:almostmaximal} is an optimal graph with 6 vertices and 17 arcs, whose algebraic connectivity is $3.215$, while the algebraic connectivity of the graph constructed by Algorithm~1, $\bbb G(6,17)$, is 3. 
In cases when the maximum algebraic connectivity is small, such a gap may be considered relatively large, and the constructed algebraic connectivity is ``far from'' the maximum. For instance, $\bbb G(4,5)$, $\bbb G(5,7)$, and $\bbb G(6,9)$ all have an algebraic connectivity of $1$, while the maximum algebraic connectivities identified through exhaustive simulation search in all three cases are $1.5$.


\begin{figure}[!ht]
\centering
\includegraphics[width=2.2in]{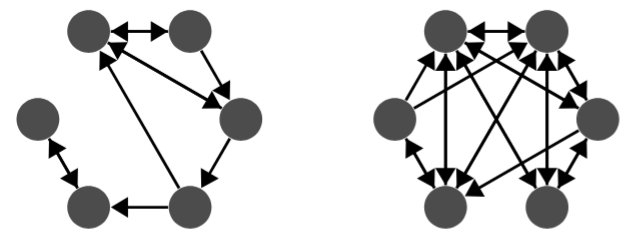} 
\caption{Two almost regular graphs with algebraic connectivity far from the maximum}
\label{fig:almostregular-smallaG}
\end{figure}

It is worth emphasizing that the property in Theorem \ref{th:p1-construction} does not hold for all almost regular graphs. To illustrate this, we provide two examples of almost regular graphs whose algebraic connectivity is far from the maximum possible value, with a gap exceeding $1$. 
The left graph in Figure \ref{fig:almostregular-smallaG} is an almost regular graph with $6$ vertices and $10$ arcs, having an algebraic connectivity of $0.161$. In contrast, the maximum algebraic connectivity in this case is $2$ and can be achieved by $\bbb G(6,10)$, as explained immediately following the statement of Theorem \ref{th:p1-construction}.
Similarly, the right graph in Figure \ref{fig:almostregular-smallaG} is an almost regular graph with $6$ vertices and $21$ arcs, with an algebraic connectivity of $2.382$, whereas the maximum algebraic connectivity of $4$ can be achieved by $\bbb G(6,21)$, as identified through an exhaustive simulation search. Both cases admit multiple optimal graphs.

To prove Theorem \ref{th:p1-construction}, we need the following lemmas regarding the relationship between $\kappa=\lfloor \frac{m}{n-1} \rfloor$ and $\nu=\lfloor \frac{m}{n} \rfloor$. 

\vspace{.05in}

\begin{lemma}\label{lm:p-k}
$\nu \in \{\kappa-1,\kappa\}$ for any integers $n$ and $m$ such that $n \ge 2$ and $1\le m \le n(n-1)$.
\end{lemma}

\vspace{.05in}

{\bf Proof of Lemma \ref{lm:p-k}:} Since $\lfloor \frac{m}{n-1} \rfloor \le \frac{m}{n-1}$ and $\lfloor \frac{m}{n} \rfloor > \frac{m}{n}-1$, it follows that 
$\kappa-\nu = \lfloor \frac{m}{n-1} \rfloor - \lfloor \frac{m}{n} \rfloor < \frac{m}{n-1} - (\frac{m}{n}-1 ) = \frac{m}{n(n-1)} + 1 \le 2$. 
As $\kappa-\nu$ is a nonnegative integer, it can only take a value of either 0 or 1, which implies that $\nu$ is equal to either $\kappa$ or $\kappa-1$.
\hfill $\qed$

\vspace{0.05in}

\begin{lemma}\label{lm:floor}
$\lfloor \frac{m-\nu-1}{n-2} \rfloor = \kappa$ for any integers $n$ and $m$ such that 
$n \ge 3$ and $1 \le m \le n(n-2)$.
\end{lemma}


{\bf Proof of Lemma \ref{lm:floor}:}
First, consider the special case when $1 \le m \le n-2$, which implies $\kappa = \nu = 0$. Then, $\lfloor \frac{m-\nu-1}{n-2} \rfloor = \lfloor \frac{m-1}{n-2} \rfloor = 0 = \kappa$. Thus, the lemma holds in this case. 

Next, consider the general case when 
$n-1 \le m \le n(n-2)$.
Note that $m$ can be written as $m = \kappa(n-1)+r$, where $0 \le r \le n-2$ is the unique remainder when $m$ is divided by $n-1$. 
From Lemma \ref{lm:p-k}, $\nu$ equals either $\kappa$ or $\kappa-1$. Let us first suppose $\nu=\kappa$. Then, $\frac{\kappa (n-1) + r}{n} = \frac{m}{n} \ge \lfloor \frac{m}{n} \rfloor =\nu =\kappa$, which implies $r \ge \kappa = \lfloor \frac{m}{n-1} \rfloor \ge 1$. Thus, 
$\lfloor \frac{m-\nu-1}{n-2} \rfloor = \lfloor \frac{\kappa(n-1)+r-\kappa-1}{n-2} \rfloor = \kappa + \lfloor \frac{r-1}{n-2} \rfloor = \kappa$. 
In the next step, we suppose $\nu=\kappa-1$. Then, $\lfloor \frac{m-\nu-1}{n-2} \rfloor = \lfloor \frac{\kappa(n-1)+r-\kappa}{n-2} \rfloor = \kappa + \lfloor \frac{r}{n-2} \rfloor$, which equals $\kappa$ if $0\le r< n-2$. 
To complete the proof, it remains to consider the case when $r=n-2$. We claim that $r\neq n-2$. To prove the claim, suppose to the contrary that $r = n-2$, with which $m = \kappa(n-1)+r = \kappa n+(n-2-\kappa)$.
Meanwhile, 
as $n \ge 3$ and $1 \le m \le n(n-2)$,
$\kappa = \lfloor \frac{m}{n-1} \rfloor \le \frac{n(n-2)}{n-1}<n-1$, which implies that $n-1-\kappa$ is a positive integer. Then, 
$n-2-\kappa \ge 0$, and thus $\nu = \lfloor \frac{m}{n} \rfloor = \lfloor \frac{\kappa n+(n-2-\kappa)}{n} \rfloor = \kappa+\lfloor \frac{n-2-\kappa}{n} \rfloor = \kappa$. 
But this contradicts $\nu=\kappa-1$. Therefore, $r\neq n-2$. 
\hfill $\qed$

\vspace{.05in}

\begin{lemma}\label{lm:n-divides-m}
For any integers $n$ and $m$ such that 
$n \ge 3$ and $1 \le m < n(n-2)$,
if $n$ divides $m$, then $\kappa=\nu=\lfloor \frac{m-\nu}{n-2} \rfloor$.
\end{lemma}

\vspace{.05in}

{\bf Proof of Lemma \ref{lm:n-divides-m}:}
Since $n$ divides $m$, $m=\nu n$. 
Then, $\kappa = \lfloor \frac{m}{n-1} \rfloor = \lfloor \frac{\nu n}{n-1} \rfloor = \lfloor \nu + \frac{\nu}{n-1} \rfloor = \nu + \lfloor  \frac{\nu}{n-1} \rfloor =\nu$. Note that $\nu = \frac{m}{n} < n-2$. Therefore,
$\lfloor \frac{m-\nu}{n-2} \rfloor = \lfloor \frac{n\nu-\nu}{n-2} \rfloor = \lfloor \nu + \frac{\nu}{n-2} \rfloor = \nu+ \lfloor  \frac{\nu}{n-2} \rfloor = \nu = \kappa$.
\hfill $\qed$

\vspace{.05in}

We will also need the following entry-wise property of the Laplacian matrix of $\bbb G(n,m)$. 
Let $L(n,m)$ denote the Laplacian matrix of $\bbb G(n,m)$, with its $ij$th entry denoted by $[L(n,m)]_{ij}$. 

\vspace{.05in}





\begin{lemma}\label{lem:Lentries}
For any integers $n$ and $m$ such that 
$n \ge 2$ and $1 \le m \le (n-1)^2$, $[L(n,m)]_{in}=0$ for each $i\in\{1,\ldots,n-1\}$.
\end{lemma}

\vspace{.05in}

{\bf Proof of Lemma \ref{lem:Lentries}:}
To prove the lemma, suppose to the contrary that $[L(n,m)]_{in} \neq 0$ for some $i\in\{1,\ldots,n-1\}$. From the definition of a Laplacian matrix, the off-diagonal entry $[L(n,m)]_{in} = -1$, which implies that $(n,i)$ is an arc in $\bbb G(n,m)$. From \eqref{eq:allGarcs}, $n = \lceil \frac{j}{n-1} \rceil$ for some $j \in \{1,\ldots,m\}$.
Then, $n = \lceil \frac{j}{n-1} \rceil \le \lceil \frac{m}{n-1} \rceil < \frac{m}{n-1} + 1$, which implies $m > (n-1)^2$. But this contradicts $m \le (n-1)^2$. Therefore, $[L(n,m)]_{in}=0$ for all $i\in\{1,\ldots,n-1\}$.
\hfill$\qed$

\vspace{.05in}

Lemma \ref{lem:Lentries} has the following implication. 

\vspace{.05in}

\begin{lemma}\label{lm:delete-vertex-n}
For any integers $n$ and $m$ such that $n \ge 3$ and $1 \le m \le (n-1)^2$, $\bbb G(n-1, m-d_n)$ is the subgraph of $\bbb G(n,m)$ induced by the vertex subset $\{1,\ldots,n-1\}$, where $d_n$ is the in-degree of vertex $n$ in $\bbb G(n,m)$. 
\end{lemma}





\vspace{.05in}

It is clear that $m-d_n\ge 0$. In the following proof, we will soon show that $m-d_n\le (n-1)(n-2)$. These two facts guarantee that $\bbb G(n-1,m-d_n)$ is well-defined.

\vspace{.05in}

\textbf{Proof of Lemma \ref{lm:delete-vertex-n}:}
Let $\bbb H$ be the subgraph of $\bbb G(n,m)$ induced by the vertex subset $\{1,\ldots,n-1\}$. From Lemma~\ref{lem:Lentries}, vertex $n$ has no outgoing arcs in $\bbb G(n,m)$. Then, $\bbb H$ has $n-1$ vertices and $m-d_n$ arcs, which implies $m-d_n\le (n-1)(n-2)$.  
In addition, each vertex $i\in\{1,
\ldots,n-1\}$ has the same in-degree $d_i$ in $\bbb H$ as in $\bbb G(n,m)$, with the values of $d_i$, $i\in\{1,\ldots,n\}$ being given in \eqref{eq:indegree}. Let $b_i$ denote the in-degree of vertex $i$ in $\bbb G(n-1,m-d_n)$. 
From Proposition~\ref{prop:p1-deg}, 
$$(b_1, \ldots, b_{n-1}) = (\underbrace{u, \ldots, u,}_{(n-1)(u+1)-(m-d_n)} \underbrace{u+1, \ldots, u+1}_{(m-d_n) - (n-1)u} \;), $$
where $u = \lfloor \frac{m-d_n}{n-1} \!\rfloor$. We claim that $b_i=d_i$ for all $i\in\{1,\ldots,n-1\}$. To prove the claim, we consider two scenarios separately.
First, suppose that $n$ divides $m$. Then, from \eqref{eq:indegree}, all $d_i$, $i\in\{1,\ldots,n\}$ equal $\nu=\frac{m}{n}$, which implies that $n-1$ divides $m-d_n$ and $u=\nu$. It follows that all $b_i$, $i\in\{1,\ldots,n-1\}$ equal $u$, and thus the claim holds. 
Next, suppose that $n$ does not divide $m$. Then, $d_n=\nu +1=\lfloor\frac{m}{n}\rfloor+1$ and $m=\nu n+r$, where $1\le r\le n-1$ is the unique remainder when $m$ is divided by $n$. With these, $u = \lfloor \frac{m-d_n}{n-1} \!\rfloor = \lfloor \frac{\nu n+r-\nu -1}{n-1}\rfloor = \nu + \lfloor \frac{r-1}{n-1} \rfloor = \nu$ and thus $(n-1)(u+1)-(m-d_n)=n(\nu+1)-m$, which validates the claim. The claim ensures that each vertex has the same number of incoming arcs in $\bbb H$ and $\bbb G(n-1,m-d_n)$. Lemma~\ref{lem:arctail} further guarantees that all these arcs coincide. Therefore, $\bbb H = \bbb G(n-1,m-d_n)$. 
\hfill $\qed$


\vspace{.05in}

We are in a position to prove Theorem \ref{th:p1-construction}, in which the following well-known result will be used. 

\vspace{.05in}

\begin{lemma}\label{lem:completegraph}
(Lemma 2 in \cite{Anderson85})
The Laplacian matrix of the $n$-vertex complete graph has a single eigenvalue of $0$ and an eigenvalue of $n$ with multiplicity $n-1$.
\end{lemma}


\vspace{.05in}

{\bf Proof of Theorem \ref{th:p1-construction}:}
It is easy to see that $m-\kappa(n-1)$ is a nonnegative integer and $(\kappa+1)(n-1)-m$ is a positive integer. 
In the special case when $m=0$, $\bbb G(n,m)$ is an empty graph and $\lfloor \frac{m}{n-1} \rfloor =0$; then the theorem is clearly true. 
For the general case when $m>0$, we will establish the following claim. 

{\bf Claim:} For any $n\ge 2$ and $1\le m\le n(n-1)$,
$$p_{L(n,m)}(\lambda) = \lambda(\lambda-\kappa)^{(\kappa+1)(n-1)-m} (\lambda-\kappa-1)^{m-\kappa(n-1)}.$$ 

Note that $1+[(\kappa+1)(n-1)-m] + [m-\kappa(n-1)] = n$. The claim implies that $L(n,m)$ has one eigenvalue at 0, $(\kappa+1)(n-1)-m$ eigenvalues at $\kappa$, and $(\kappa+1)(n-1)-m$ at $\kappa+1$, which together constitute the entire spectrum of $L(n,m)$. The theorem then immediately follows from the claim. Thus, to prove the theorem, it is sufficient to establish the claim. 
We will prove the claim by induction on $n$.


In the base case when $n=2$, all possible values of $m$ are $1$ and $2$. According to the algorithm description, $\bbb G(2,1)$ contains one arc, $(1,2)$, and $\bbb G(2,2)$ contains two arcs, $(1,2)$ and $(2,1)$. It is straightforward to verify that the claim holds for both $L(2,1)$ and $L(2,2)$.  

For the inductive step, suppose that the claim holds for $n=q$ with all possible values of $m$ in $\{1,\ldots,q(q-1)\}$, where $q\ge 2$ is an integer. 
Let $n=q+1$, and thus all possible values of $m$ range from 1 to $(q+1)q$.


We first consider the case when $m\in\{1,\ldots,q^2-1\}$.
From Lemma \ref{lem:Lentries}, with $n$ replaced by $q+1$, $[L(q+1,m)]_{i(q+1)}=0$ for all $i\in\{1,\ldots,q\}$. That is, all entries in the $(q+1)$th column of $L(q+1,m)$, except for the last entry $[L(q+1,m)]_{(q+1)(q+1)}$, are zero. The same holds for the $(q+1)$th column of $\lambda I-L(q+1,m)$, whose last entry is equal to $\lambda - [L(q+1,m)]_{(q+1)(q+1)}$. 
Then, the Laplace expansion along the $n$th column of $L(q+1,m)$ yields
\begin{align}
p_{L(q+1,m)}(\lambda) &= \det\big(\lambda I - L(q+1,m)\big) \nonumber\\
&= \big(\lambda - [L(q+1,m)]_{(q+1)(q+1)}\big) p_M(\lambda),\label{eq:characteristic}
\end{align}
where $M$ is the $q\times q$ submatrix of $L(q+1,m)$ obtained by removing the $(q+1)$th row and $(q+1)$th column of $L(q+1,m)$. 
Since all row sums of $L(q+1,m)$ are zero, and all entries in its $(q+1)$th column, except for the last entry, are zero, $M$ also has all row sums equal to zero. It follows that $M$ is the Laplacian matrix of a certain graph $\bbb H$ with $q$ vertices, and $\bbb H$ is the subgraph of $\bbb G(q+1,m)$ induced by the vertex subset $\{1,\ldots,q\}$. 
From Lemma \ref{lm:delete-vertex-n}, with $n$ replaced by $q+1$, $\bbb H=\bbb G(q, m-d_{q+1})$, and thus $M$ is the Laplacian matrix of $\bbb G(q, m-d_{q+1})$, where $d_{q+1}$ denotes the in-degree of vertex $q+1$ in $\bbb G(q+1,m)$. 
Let $\gamma \dfb \lfloor \frac{m}{q+1} \rfloor$. We consider the following two cases separately.



{\bf Case 1:} Suppose that $q+1$ divides $m$, which implies $m = \gamma (q+1)$. From the definition of a Laplacian matrix and Proposition \ref{prop:p1-deg}, $[L(q+1,m)]_{(q+1)(q+1)} = d_{q+1} = \gamma$. Then, $M$ is the Laplacian matrix of $\bbb G(q, m-d_{q+1})=\bbb G(q,m-\gamma)$. From \eqref{eq:characteristic} and the induction hypothesis,  
\eq{
p_{L(q+1,m)}(\lambda) 
= (\lambda - \gamma) p_M(\lambda),\label{eq:case1a}
}
$$p_M(\lambda) = \lambda(\lambda-\beta)^{(\beta+1)(q-1)-m+\gamma} (\lambda-\beta-1)^{m-\gamma-\beta(q-1)},$$
where $\beta \dfb \lfloor\frac{m-\gamma}{q-1}\rfloor$.
The analysis is further divided into two scenarios based on the value of $m$. First, consider when $m = q^2-1$, which implies $\gamma=q-1$ and thus $\beta=q$. Then, $(\beta+1)(q-1)-m+\gamma = q-1$ and $m-\gamma-\beta(q-1) = 0$. 
It follows from \eqref{eq:case1a} that
\eq{
p_{L(q+1,m)}(\lambda) = \lambda(\lambda-q+1)(\lambda -1)^{q-1}. \nonumber 
}
Meanwhile, $(\gamma+1)q-m = 1$ and $m-\gamma q = q-1$. Thus, the above equation validates the claim with $n$ replaced by $q+1$. 
Next, consider when $1\le m \le q^2-2$. From Lemma~\ref{lm:n-divides-m}, with $n$ replaced by $q+1$, $\gamma = \beta =\alpha \dfb \lfloor\frac{m}{q}\rfloor$. Then, from~\eqref{eq:case1a},
\eq{
p_{L(q+1,m)}(\lambda) = \lambda(\lambda-\alpha)^{(\alpha+1)q-m}(\lambda -\alpha-1)^{m-\alpha q},\label{eq:consistentclaim}
}
which proves the claim with $n$ replaced by $q+1$.

{\bf Case 2:} Suppose that $q+1$ does not divide $m$, which implies $m -\gamma (q+1)>0$. 
From the definition of a Laplacian matrix and Proposition \ref{prop:p1-deg}, $[L(q+1,m)]_{(q+1)(q+1)} = d_{q+1} = \gamma +1$. Then, $M$ is the Laplacian matrix of $\bbb G(q, m-d_{q+1})=\bbb G(q,m-\gamma-1)$. From \eqref{eq:characteristic} and the induction hypothesis,  
\eq{
p_{L(q+1,m)}(\lambda) 
= (\lambda - \gamma-1) p_M(\lambda),\label{eq:case1b}
}
$$p_M(\lambda) \! = \! \lambda(\lambda-\beta')^{(\beta'+1)(q-1)-m+\gamma'} (\lambda-\beta'-1)^{m-\gamma'-\beta'(q-1)},$$
where $\gamma'=\gamma+1$ and $\beta' \dfb \lfloor\frac{m-\gamma-1}{q-1}\rfloor$.
With $n$ replaced by $q+1$, Lemma \ref{lm:floor} and Lemma \ref{lm:p-k} respectively imply that $\beta'=\alpha$ and $\gamma\in\{\alpha-1,\alpha\}$. 
The analysis is then divided into two scenarios based on the value of $\gamma$. 
First, suppose $\gamma = \alpha-1$. Then, $(\beta'+1)(q-1)-m+\gamma' = (\alpha+1)(q-1)-m+\alpha$ and $m-\gamma'-\beta'(q-1) = m-\alpha q$. It follows that \eqref{eq:case1b} simplifies to \eqref{eq:consistentclaim}, 
which validates the claim.  
Next, suppose $\gamma = \alpha$. 
Then, $(\beta'+1)(q-1)-m+\gamma' = (\alpha+1)q-m$ and $m-\gamma'-\beta'(q-1) = m-\alpha q-1$. With these equalities, \eqref{eq:case1b} once again leads to \eqref{eq:consistentclaim}, thereby proving the claim.

The two cases above collectively establish the inductive step for $m\in\{1,\ldots,q^2-1\}$.
In what follows, we address the scenario where $m \in\{q^2,\ldots, q^2+q\}$.
In the special case when $m=q^2+q$, $\bbb G(q+1,m)$ is the complete graph. From Lemma \ref{lem:completegraph}, $\bbb G(q+1,m)$ has a single eigenvalue at $0$ and $q$ eigenvalues at $q+1$, which implies $p_{L(q+1,m)}(\lambda) = \lambda (\lambda-q-1)^{q}$. Note that in this case $\alpha = q+1$ and $m-\alpha q = 0$. The characteristic polynomial expression proves the claim.

It remains to consider when $m \in \{q^2,\ldots,q^2+q-1\}$, with which $\alpha =q$. From Proposition \ref{prop:p1-stars}, with $n$ replaced by $q+1$, $\bbb G(q+1,m)$ is the union of $q$ simple $(q+1)$-vertex directed stars, rooted at vertices $1, \ldots, q$, respectively, and the $(q+1)$-vertex simple directed forest consisting of an $(m-q^2+1)$-vertex directed star, rooted at vertex $q+1$, and $q^2+q-m$ isolated vertices.
Then, $\overline{\bbb G(q+1,m)}$ is a $(q+1)$-vertex simple directed forest consisting of a $(q^2+q-m+1)$-vertex directed star, rooted at vertex $q+1$, and $m-q^2$ isolated vertices. Thus, $\overline{\bbb G(q+1,m)}$ is acyclic, containing $q^2+q-m$ vertices with an in-degree of $1$, while all other vertices have an in-degree of $0$.
From Lemma \ref{lm:acyclic}, the Laplacian spectrum of $\overline{\bbb G(q+1,m)}$ consists of $q^2+q-m$ eigenvalues at $1$ and $m-q^2+1$ eigenvalues at $0$. Then, from Lemma \ref{lm:complement}, the Laplacian spectrum of $\bbb G(q+1,m)$ is composed of a single eigenvalue of $0$, an eigenvalue of $q$ with multiplicity $q^2+q-m$, and an eigenvalue of $q+1$ with multiplicity $m-q^2$. This leads to the characteristic polynomial 
$$p_{L(q+1,m)}(\lambda) = \lambda (\lambda-q)^{q^2+q-m}(\lambda-q-1)^{m-q^2},$$
which validates the claim with $n$ replaced by $q+1$. This completes the proof of the inductive step.
\hfill $\qed$

\section{Conclusion}

In this paper, we have studied a couple of special cases of an open problem: identifying optimal directed graphs for the fastest convergence rate of continuous-time consensus. Necessary and sufficient conditions for both sparse and dense graphs to be optimal have been derived. The optimal consensus performance of any union of directed stars (with the same vertex set) has also been established.

For general cases involving arbitrary numbers of vertices and arcs, the open problem remains quite challenging. To provide an approximate solution to the problem, a computationally efficient algorithm has been proposed. The algorithm employs an inductive construction process to generate a sequence of almost regular directed graphs for each possible number of arcs, given the number of vertices. These graphs attain algebraic connectivity that is close to the maximum possible value, with a gap of at most one. Consequently, the constructed almost regular directed graphs guarantee fast consensus performance for continuous-time consensus processes implemented on them.
The proposed algorithm can thus be regarded as a generalization of the algorithm presented in \cite{cdc24laplacian}, which generates almost regular (undirected) graphs for any specified pair of vertex and edge numbers, achieving large algebraic connectivity.

One direction for future work is to investigate vertex and (directed) edge connectivity in the almost regular directed graphs constructed by the proposed algorithm, similar to what has been done for undirected almost regular graphs in \cite{cdc24laplacian}.
Another important direction is to apply or extend the tools developed here and in \cite{cdc24laplacian} from continuous-time to discrete-time linear consensus processes \cite{vicsekmodel}, as identifying the fastest consensus graphs for the latter, both undirected and directed, also remains largely unexplored.

\begin{figure}[!ht]
\centering
\includegraphics[width=2.2in]{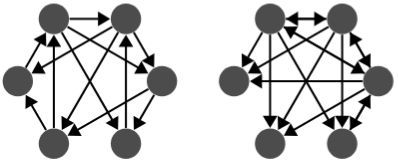} 
\caption{Two optimal graphs with maximal algebraic connectivity}
\label{fig:fastest-consensus-ex}
\end{figure}

Undirected graphs can be considered a special subclass of directed graphs, with each undirected edge equivalently represented by a pair of arcs with opposite directions; hence, they are also referred to as symmetric directed graphs. The results presented in this paper partially demonstrate that, given a fixed number of vertices and arcs, the optimal graphs corresponding to the fastest consensus are almost never symmetric/undirected. We provide two concrete examples in Figure \ref{fig:fastest-consensus-ex} to illustrate this point. Both graphs in the figure have 6 vertices. All $6$-vertex optimal graphs with the largest algebraic connectivity for different numbers of undirected edges are illustrated in \cite[Figure 4]{ogiwara}. The graph on the left achieves the maximal algebraic connectivity of $2.215$ among all directed graphs with $6$ vertices and $12$ arcs; in contrast, the maximum algebraic connectivity among all undirected graphs with $6$ vertices and $6$ undirected edges (equivalent to $12$ arcs) is $1$. Similarly, the graph on the right, which is $\bbb G(6,16)$ constructed by Algorithm 1, achieves the maximal algebraic connectivity of $3$ among all directed graphs with $6$ vertices and $16$ arcs, whereas the maximum algebraic connectivity among all undirected graphs with $6$ vertices and $8$ undirected edges (equivalent to $16$ arcs) is $2$. These observations may call into question the view that consensus is generally faster in undirected graphs compared to directed~graphs.

Finally, the results in this paper are closely related to an unsolved conjecture in network synchronization.
Let $L$ be the Laplacian matrix of an $n$-vertex simple directed graph $\bbb G$, whose eigenvalues are $\lambda_1,\lambda_2,\ldots,\lambda_n$ with $\lambda_1=0$. 
Define
$$\textstyle\sigma^2 \dfb \frac{1}{n-1} \sum_{i=2}^n |\lambda_i - \bar\lambda|^2, \;\;\;\; {\rm where} \;\;\;\; \bar\lambda = \frac{1}{n-1} \sum_{i=2}^n \lambda_i,$$
which is a normalized deviation of possibly nonzero eigenvalues. This quantity is called the normalized spread of the eigenvalues in \cite{Nishikawa10} to measure the synchronizability of certain network dynamics. It is claimed and validated by simulations that the smaller the value of $\sigma^2n^2/m^2$, the more synchronizable the network will generally be, where $m$ denotes the number of directed edges in $\bbb G$.  
It is conjectured in \cite{Nishikawa10}, without theoretical validation, that 
among all simple directed graphs with $n$ vertices and $m$ arcs, the minimum possible value of $\sigma^2$ is
achieved if the Laplacian spectrum is 
$$
0, \underbrace{\kappa,\; \ldots, \;\kappa,}_{(\kappa+1)(n-1)-m} \; \underbrace{\kappa+1, \ldots, \kappa+1}_{m-\kappa(n-1)}.
$$
We have shown that the graph constructed by Algorithm 1, $\bbb G(n,m)$, has the same Laplacian spectrum as conjectured above. Therefore, this paper has proved that the conjectured minimum value of $\sigma^2$ is always achievable for any feasible pair of $n$ and $m$. The question of whether this value of $\sigma^2$ is truly the minimum remains an open problem.


\section*{Acknowledgement}
The authors wish to thank Yuanyuan Yang (Stony Brook University) and Yiming Zeng (Binghamton University) for their valuable discussions and efforts in searching for large-scale optimal graphs via cloud computing.

\bibliographystyle{unsrt}
\bibliography{susie,consensus,jicareer}

\end{document}